\documentclass[12pt,draftcls,onecolumn]{IEEEtran}

\usepackage{amsmath,amsfonts,amssymb}
\usepackage{tikz, subcaption}
\newtheorem{definition}{Definition}
\newtheorem{theorem}{Theorem}
\newtheorem{lemma}{Lemma}

\newtheorem{corollary}{Corollary}
\newtheorem{remark}{Remark}
\newtheorem{example}{Example}

\usepackage{epsfig} % for postscript graphics files

\title{On the Smallest Eigenvalue of\\ Grounded Laplacian Matrices}

\author{Mohammad~Pirani~and~Shreyas~Sundaram% <-this % stops a space
\thanks{This material is based upon work supported in part by the Natural Sciences and Engineering Research Council of Canada (NSERC).  The authors are with the Department of Electrical and Computer Engineering, University of Waterloo, Canada. E-mail for corresponding author: ssundara@uwaterloo.ca. }%
}

\begin{document}

\maketitle
\thispagestyle{empty}
\pagestyle{empty}

\begin{abstract}
We provide bounds on the smallest eigenvalue of grounded Laplacian matrices (which are obtained by removing certain rows and columns of the Laplacian matrix of a given graph).  The gap between our upper and lower bounds depends on the ratio of the smallest and largest components of the eigenvector corresponding to the smallest eigenvalue of the grounded Laplacian.   We provide a graph-theoretic bound on this ratio, and subsequently obtain a tight characterization of the smallest eigenvalue for certain classes of graphs.  Specifically,  for Erdos-Renyi random graphs, we show that when a (sufficiently small) set $S$ of rows and columns is removed from the Laplacian, and the probability $p$ of adding an edge is sufficiently large, the smallest eigenvalue of the grounded Laplacian asymptotically almost surely approaches $|S|p$.  We also show that for random $d$-regular graphs with a single row and column removed, the smallest eigenvalue is $\Theta(\frac{d}{n})$.  Our bounds have applications to the study of the convergence rate in consensus dynamics with stubborn or leader nodes.
\end{abstract}

\iffalse
\fi

\section{Introduction}
 \label{sec:intro}

There has been a great deal of research over the past several decades dedicated to the study of the structure and dynamics of networks.  These investigations span multiple disciplines and include combinatorial, probabilistic, game-theoretic, and algebraic perspectives \cite{Harary, Bollobas, JacksonBook, Godsil}.  It has been recognized that the {\it spectra} of graphs (i.e., the eigenstructure associated with certain matrix representations of the network) provide insights into both the topological properties of the underlying network and dynamical processes occurring on the network \cite{ChungSpectral, Cvetkovic}. The eigenvalues and eigenvectors of the Laplacian matrix of the graph, for example, contain information about the connectivity and community structure of the network \cite{Fiedler,Merris,Mohar,Luxburg}, and dictate the convergence properties of certain diffusion dynamics \cite{Olfati1}.
 
A variant on the Laplacian that has attracted attention in recent years is the {\it grounded Laplacian} matrix, obtained by removing certain rows and columns from the Laplacian.  The grounded Laplacian forms the basis for the classical Matrix Tree Theorem (characterizing the number of spanning trees in the graph), and also plays a fundamental role in the study of continuous-time diffusion dynamics where the states of some of the nodes in the network are fixed at certain values.   The eigenvalues of the grounded Laplacian characterize the variance in the equilibrium values for noisy instances of such dynamics, and determine the rate of convergence to steady state \cite{Bamieh,Clark1}.  Optimization algorithms have been developed to select ``leader nodes'' in the network in order to minimize the steady-state variance or to maximize the rate of convergence \cite{Clark1,Clark2,Jovanovic,Leonard}, and various works have studied the effects of the location of such leaders in distributed control and consensus dynamics \cite{Barooah,Hao,Shi13,Ghaderi13,ACC}.  
 
In this paper, we provide a characterization of the smallest eigenvalue of grounded Laplacian matrices.  Specifically, we provide graph-theoretic bounds on the smallest eigenvalue based on the number of edges leaving the grounded nodes, bottlenecks in the graph, and properties of the eigenvector associated with the eigenvalue.  Our bounds become tighter as this eigenvector becomes more uniform; we provide graph properties under which this occurs.  As a consequence of our analysis, we obtain the smallest eigenvalue of the grounded Laplacian matrix for Erdos-Renyi random graphs and random regular graphs.

\section{Background and Notation}
\label{sec:defs}
We use $\mathcal{G} = \{\mathcal{V},\mathcal{E}\}$ to denote an undirected graph where $\mathcal{V}$ is the set of nodes (or vertices) and $\mathcal{E} \subset \mathcal{V}\times\mathcal{V}$ is the set of edges.  We will denote the number of vertices by $n$.  The {\it neighbors} of  node $v_i \in \mathcal{V}$ in graph $\mathcal{G}$ are given by the set $\mathcal{N}_i = \{v_j \in \mathcal{V}~|~(v_i, v_j) \in \mathcal{E}\}$.  The {\it degree} of node $v_i$ is $d_i = |\mathcal{N}_i|$, and the minimum and maximum degrees of the nodes in the graph will be denoted by $d_{min}$ and $d_{max}$, respectively.  If $d_{max} = d_{min} = d$, the graph is said to be $d$-{\it regular}.  For a given set of nodes $S \subset \mathcal{V}$, the {\it edge-boundary} (or just boundary) of the set is given by $\partial{S} = \{(v_i,v_j) \in \mathcal{E} \mid v_i \in S, v_j \in \mathcal{V}\setminus{S}\}$.  The {\it isoperimetric constant} of $\mathcal{G}$ is given by \cite{ChungSpectral}
 $$
 i(\mathcal{G})\triangleq \min_{A \subset \mathcal{V}, |A| \le \frac{n}{2}}\frac{|\partial A|}{|A|}.
 $$
 Choosing $A$ to be the vertex with the smallest degree yields the bound $i(\mathcal{G}) \le d_{min}$.
 
 \subsection{Laplacian and Grounded Laplacian Matrices}
 \label{subsec:prop_grounded_laplacian}
 The {\it adjacency matrix} for the graph is a matrix $A \in \{0,1\}^{n\times{n}}$, where entry $(i,j)$ is $1$ if $(v_i,v_j) \in \mathcal{E}$, and zero otherwise.  The {\it Laplacian matrix} for the graph is given by $L = D - A$, where $D$ is the degree matrix with $D = \hbox{diag}(d_1, d_2, \ldots, d_n)$.  For an undirected graph $\mathcal{G}$, the Laplacian $L$ is a symmetric matrix with real eigenvalues that can be ordered sequentially as $0=\lambda_1 (L)\leq \lambda_2 (L)\leq \cdots \leq \lambda_n (L)\leq 2d_{max}$.  The second smallest eigenvalue $\lambda_2(L)$ is termed the {\it algebraic connectivity} of the graph and satisfies the bound \cite{ChungSpectral} 
 \begin{equation}
\lambda_2(L) \ge \frac{i(\mathcal{G})^2}{2d_{max}}.
\label{eqn:lower_bound_lambda_2_iso}
\end{equation}

  We will designate a nonempty subset of  vertices $\mathcal{S} \subset \mathcal{V}$ to be {\it grounded nodes}, and assume without loss of generality that they are placed last in the ordering of the nodes. 
 We use $\alpha_i$ to denote the number of grounded nodes that node $v_i$ is connected to (i.e., $\alpha_i = \left|\mathcal{N}_i \cap \mathcal{S}\right|$).  
 Removing the rows and columns of  $L$ corresponding to the grounded nodes $S$ produces a {\it grounded Laplacian matrix} (also known as a {\it Dirichlet Laplacian matrix}) denoted by $L_{g}(S)$.  When the set $S$ is fixed and clear from the context, we will simply use $L_g$ to denote the grounded Laplacian.
For any given set $S$, we denote the smallest eigenvalue of the grounded Laplacian by $\lambda(L_{g}(S))$ or simply $\lambda$.  

When the graph $\mathcal{G}$ is connected, the grounded Laplacian matrix is a positive definite matrix and all of the elements in its inverse are nonnegative \cite{Miekkala}.  From the Perron-Frobenius (P-F) theorem \cite{Horn}, the eigenvector associated with the smallest eigenvalue of the grounded Laplacian can be chosen to be nonnegative (elementwise).  Furthermore, when the grounded nodes do not form a vertex cut, the eigenvector associated with the smallest eigenvalue is unique (up to normalization) can be chosen to have all elements positive.

%%%%%%%%%%%%%%%%%%%%%%%%%%%%%%%%%%%%%%%%%%%%%%%%%%%%%%%%%%%%%%%%%%%
\subsection{Applications in Consensus with Stubborn Agents} 
Consider a multi-agent system described by the connected and undirected graph $\mathcal{G} = \{\mathcal{V},\mathcal{E}\}$ representing the structure of the system, and a set of equations describing the interactions between each pair of agents.  In the study of consensus and opinion dynamics \cite{Olfati1}, each agent $v_i \in \mathcal{V}$ starts with an initial scalar state (or opinion) $y_i(t)$, which evolves over time as a function of the states of its neighbors.  A commonly studied version of these dynamics involves a continuous-time linear update rule of the form
\begin{equation*}
\dot{y}_i(t)= \sum_{v_j\in \mathcal{N}_i} (y_j(t)-y_i(t)).
%\label{eqn:dynamic}
\end{equation*}
Aggregating the state of all of the nodes into the vector $Y(t)=\left[\begin{matrix}y_1(t) & y_2(t) & \cdots & y_n(t)\end{matrix}\right]^T$,  the above equation produces the system-wide dynamical equation
\begin{equation}
\dot{Y}=-LY,
\label{eqn:dynamicmatrix}
\end{equation}
where $L$ is the graph Laplacian.   When the graph is connected, the trajectory of the above dynamical system satisfies $Y(t) \rightarrow \frac{1}{n}\mathbf{1}\mathbf{1}^TY(0)$ (i.e., all agents reach consensus on the average of the initial values), and the asymptotic convergence rate is given by $\lambda_2(L)$ \cite{Olfati1}.

Now suppose that there is a subset $S \subset \mathcal{V}$ of agents whose opinions are kept constant throughout time, i.e., $\forall v_s \in S$, $\exists y_s \in \mathbb{R}$ such that $y_s(t)=y_s$ $\forall t \in \mathbb{R}_{\ge{0}}$. Such agents are known as {\it stubborn agents} or {\it leaders} (depending on the context) \cite{Clark1, Ghaderi13}. In this case the dynamics \eqref{eqn:dynamicmatrix} can be written in the matrix form
\begin{equation*}
\begin{bmatrix}
  \dot{Y}_F(t)   \\
  \dot{Y}_S(t) \\
 \end{bmatrix}=-\begin{bmatrix}
  L_{11} & L_{12}   \\
  L_{21}  & L_{22}  \\
 \end{bmatrix}\begin{bmatrix}
  {Y}_F(t)   \\
  {Y}_S(t) \\
 \end{bmatrix},
% \label{eqn:matrixstubborn}
\end{equation*}
where $Y_F$ and $Y_S$ are the states of the non-stubborn and stubborn agents, respectively. Since the stubborn agents keep their values constant, the matrices $L_{21}$ and $L_{22}$ are zero. Thus, the matrix $L_{11}$ is the grounded Laplacian for the system, i.e., $L_{11} = L_g(S)$.   It can be shown that the state of each follower asymptotically converges to a convex combination of the values of the stubborn agents, and that the rate of convergence is asymptotically given by $\lambda$, the smallest eigenvalue of the grounded Laplacian \cite{Clark1}.  

Similarly, one can consider discrete-time consensus dynamics (also known as DeGroot dynamics) with a set $S$ of stubborn nodes, given by the update equation $
Y_F(t+1)=A_gY_F(t)$, where $Y_F(t)$ is the state vector for the non-stubborn nodes at time-step $t$, and $A_g$ is an $(n-|S|)\times(n-|S|)$ nonnegative matrix given by $A_g = I - \frac{1}{k}L_g$, with constant $k \in (d_{max},\infty)$ \cite{Kingston}.
Once again, each non-stubborn node will converge asymptotically to a convex combination of the stubborn nodes' states.  The largest eigenvalue of $A_g$ is given by $\lambda_{max}(A_g)=1-\frac{1}{k}\lambda(L_g)$, and determines the asymptotic rate of convergence.  Thus, our bounds on the smallest eigenvalue of the grounded Laplacian will readily translate to bounds on the largest eigenvalue of $A_g$.

There have been various recent investigations of graph properties that impact the convergence rate for a given set of stubborn agents, leading to the development of algorithms to find approximately optimal sets of stubborn/leader agents to maximize the convergence rate \cite{Clark1, Ghaderi13, Shi13}.  
The bounds provided in this paper contribute to the understanding of consensus dynamics with fixed opinions by providing bounds on the convergence rate induced by any given set of stubborn or leader agents.

%%%%%%%%%%%%%%%%%%%%%%%%%%%%%%%%%%%%%%%%%%%%%%%%%%%%%%%%%%%%%%
\section{Bounds on the Smallest Eigenvalue of $L_g$ }
 \label{sec:spectralprop}
  
The following theorem provides our core bounds on the smallest eigenvalue of the grounded Laplacian; in subsequent sections, we will characterize graphs where these bounds become tight.

\begin{theorem}
Consider a graph $\mathcal{G}= \{\mathcal{V},\mathcal{E}\}$ with a set of grounded nodes $S \subset \mathcal{V}$.  Let $\lambda$ denote the smallest eigenvalue of the grounded Laplacian $L_g$ and let $\mathbf{x}$ be a corresponding nonnegative eigenvector, normalized so that the largest component is $x_{max} = 1$.  
Then
\begin{equation}
\frac{|\partial S|}{n-|S|}x_{min} \leq \lambda \leq \min_{X\subseteq \mathcal{V}\setminus S} \frac{|\partial X|}{|X|}\leq \frac{|\partial S |}{n-|S|},
\label{eqn:maineq}
\end{equation}
where $x_{min}$ is the smallest eigenvector component in $\mathbf{x}$.  
\label{thm:main}
\end{theorem}

\begin{IEEEproof}
From the Rayleigh quotient inequality \cite{Horn}, we have
$$
\lambda \leq z^T L_g z ,
$$
for all $z\in \mathbb{R}^{n-|S|}$ with $z^Tz=1$. Let $X \subseteq \mathcal{V}\setminus{S}$ be the subset of vertices for which $\frac{|\partial X|}{|X|}$ is minimum, and assume without loss of generality that the vertices are arranged so that those in set $X$ come first in the ordering.  The upper bound  $ \min_{X\subseteq \mathcal{V}\setminus S} \frac{|\partial X|}{|X|}$ is then obtained by choosing $z=\frac{1}{\sqrt{|X|}}[\mathbf{1}_{1\times |X|} \quad \mathbf{0}_{1\times |\mathcal{V}\setminus \{X \cup S\}|}]^T$, and noting that the sum of all elements in the top $|X| \times |X|$ block of $L_g$ is equal to the sum of the number of neighbors each vertex in $X$ has outside $X$ (i.e., $|\partial X|$).   The upper bound $\frac{|\partial S |}{n-|S|}$ readily follows by choosing the subset $X = \mathcal{V}\setminus{S}$.

For the lower bound, we left-multiply the eigenvector equation $L_g \mathbf{x} = \lambda \mathbf{x}$ by the vector consisting of all $1$'s to obtain
\begin{equation*}
\sum_{i=1}^{n-|S|}\alpha_ix_i =\lambda\sum_{i=1}^{n-|S|} x_i,
%\label{eqn:lambda}
\end{equation*}
where $\alpha_i$ is the number of grounded nodes in node $v_i$'s neighborhood.  Using the fact that the eigenvector is nonnegative, this gives
\begin{equation*}
x_{min}\sum_{i=1}^{n-|S|}\alpha_i\leq \sum_{i=1}^{n-|S|}\alpha_ix_i=\lambda\sum_{i=1}^{n-|S|} x_i \leq \lambda(n-|S|)  x_{max}=\lambda(n-|S|).
\end{equation*}
Since $\sum_{i=1}^{n-|S|}\alpha_i=|\partial S|$, the lower bound is obtained.
\end{IEEEproof}

\begin{remark}
For the case that $|S|=1$ we have 
\begin{equation*}
\frac{d_s x_{min}}{n-1}\leq \lambda\leq \frac{d_s}{n-1},
%\label{eqn:maineq2}
\end{equation*}
where $d_s$ is the degree of the grounded node.  Note that the smallest eigenvalue of the grounded Laplacian for a set $S$ of grounded nodes is always upper bounded by $|S|$ (since $|\partial S| \le |S|(n-|S|)$), with equality if and only if all grounded nodes connect to all other nodes (it is easy to see that the smallest eigenvector component $x_{min} = 1$ in this case).
\end{remark}

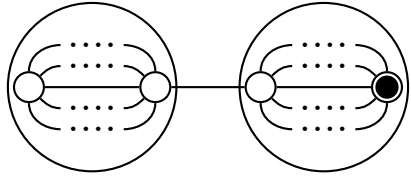
\begin{figure}[t]
\begin{center}
\begin{tikzpicture}[scale=.28, auto, node distance=1cm, thick,
   node_style/.style={circle,draw=black,fill=white!20!,font=\sffamily\small\bfseries},
   edge_style/.style={draw=black, thick}]
 \node[node_style](n1) at (9,4)  {};
 \node[node_style](n2) at (4,4)  {};
 \node[node_style](n3) at (15,4)  {};
\node[node_style](n4) at (-2,4)  {};

\draw[edge_style]  (n1) edge node{} (n2);

\draw (n4) to (n2);
\draw (n3) to (n1);

\draw    (n3) to[out=135,in=0] (13.5,5);
\draw    (n3) to[out=-135,in=0] (13.5,3);
\draw    (n3) to[out=90,in=0] (13.5,6);
\draw    (n3) to[out=-90,in=0] (13.5,2);
\draw    (n1) to[out=45,in=180] (10.5,5);
\draw    (n1) to[out=-45,in=180] (10.5,3);
\draw    (n1) to[out=-90,in=180] (10.5,2);
\draw    (n1) to[out=90,in=180] (10.5,6);
\draw    (n2) to[out=135,in=0] (2.5,5);
\draw    (n2) to[out=-135,in=0] (2.5,3);
\draw    (n2) to[out=90,in=0] (2.5,6);
\draw    (n2) to[out=-90,in=0] (2.5,2);
\draw    (n4) to[out=45,in=180] (-.5,5);
\draw    (n4) to[out=90,in=180] (-.5,6);
\draw    (n4) to[out=-90,in=180] (-.5,2);
\draw    (n4) to[out=-45,in=180] (-.5,3);

\filldraw[color=black!100, fill=black!100, thick](15,4) circle (0.5);

\node (n3) at (12,6) {{\Large $....$}};
\node (n4) at (12,5) {{\Large $....$}};
\node (n4) at (12,3) {{\Large $....$}};
\node (n4) at (12,2) {{\Large $....$}};

\node (n4) at (1,6) {{\Large $....$}};
\node (n4) at (1,5) {{\Large $....$}};
\node (n4) at (1,3) {{\Large $....$}};
\node (n4) at (1,2) {{\Large $....$}};

  \draw (1,4) ellipse (4cm and 4cm);
  \draw (12,4) ellipse (4cm and 4cm);

\end{tikzpicture}
\caption{Two complete graphs, each with $\frac{n}{2}$ nodes, connected via a single edge. The grounded node is colored black.}
\label{fig:dumbbell}
\end{center}
\end{figure}

\begin{example}
Consider the graph shown in Figure~\ref{fig:dumbbell} consisting of two complete graphs on $\frac{n}{2}$ nodes, joined by a single edge.  Suppose the black node in the figure is chosen as the grounded node.  In this case, we have $|\partial S| = \frac{n}{2} - 1$, and the extreme upper bound in \eqref{eqn:maineq} indicates that $\lambda \le \frac{|\partial S|}{n-1} \approx \frac{1}{2}$ for large $n$.  Now, if we take $X$ to be the set of all nodes in the left clique, we have $|\partial X| = 1$ and $|X| = \frac{n}{2}$, leading to $\lambda \le \frac{2}{n}$ by the intermediate upper bound in \eqref{eqn:maineq}.  
\label{ex:dumbbell}
\end{example}

In the next section, we will characterize graphs under which $x_{min}$ (the smallest eigenvector component) converges to $1$, in which case the lower and upper bounds in \eqref{eqn:maineq} coincide and yield a tight characterization of $\lambda$.  As seen in the above example, the presence of bottlenecks among the non-grounded nodes will cause $x_{min}$ to go to zero;  in certain graphs with good expansion properties, however, we will see that this will not occur.

\section{The Behavior of the Smallest Eigenvector Component}
  \label{sec:regim}
  
In this section, we analyze the effect of the network structure  on the behavior of the smallest eigenvector component $x_{min}$.  We will provide conditions under which this component goes to $1$ and stays bounded away from $0$, respectively.  This will then allow us to characterize the tightness of the bounds on the smallest eigenvalue in \eqref{eqn:maineq}.   

For a given subset $S \subset \mathcal{V}$ of grounded nodes, let $L_g(S)$ be the grounded Laplacian matrix with smallest eigenvalue $\lambda$ and corresponding nonnegative eigenvector $\mathbf{x}$.  We denote the $i$-th element of $\mathbf{x}$ by $x_i$. 
We write $L_g=\bar{L}+\Delta$ where $\bar{L}$ is the $(n-|S|)\times (n-|S|)$ Laplacian matrix of the graph when we remove the grounded nodes and all of their incident edges.  Matrix $\Delta$ is a ${(n-|S|)\times (n-|S|)}$ diagonal matrix with the $i$-th diagonal element equal to $\alpha_i$ (the number of grounded neighbors of node $v_i$).     We assume the graph corresponding to $\bar{L}$ is connected (as $x_{min}$ can be $0$ otherwise), and denote the eigenvalues of $\bar{L}$ by $0=\lambda_1(\bar{L}) <  \lambda_2(\bar{L})\leq ... \leq \lambda_{n-|S|}(\bar{L})$, with corresponding orthogonal eigenvectors $\mathbf{v}_1, \mathbf{v}_2, \ldots, \mathbf{v}_{n-|S|}$.   We take $\mathbf{v}_1 = \mathbf{1}$, and  normalize all of the other eigenvectors so that $\|\mathbf{v}_i\|_2 = 1$.  

There are various results in the literature that characterize the change in eigenvectors under modifications of matrix elements, including the commonly used Davis-Kahan theorems (which provide bounds on the angle between the original and perturbed eigenvectors) \cite{Demmel}.  However, such bounds on the angle are not particularly useful in characterizing the behavior of the smallest component of the perturbed eigenvector.\footnote{For example, consider two $n \times 1$ vectors, the first of which consists of all entries equal to $1$, and the second which has $n-1$ entries equal to $1$ and the last component equal to $0$.  The angle between these two vectors goes to $0$ as $n$ increases, but the smallest component of the second vector is always $0$.}  We thus provide the following perturbation result bounding the smallest eigenvector component of $\mathbf{x}$ in terms of the number of grounded nodes, the number of edges they have to the other nodes, and the connectivity of the graph induced by the non-grounded nodes.  The proof of the lemma starts in a similar manner to the proof of standard perturbation results \cite{Demmel}, but the latter half of the proof leverages the explicit nature of the perturbations to obtain a bound on the smallest eigenvector component (i.e., this result can be viewed as providing a bound on the $\infty$-norm of the difference between the original and perturbed eigenvectors, as opposed to a bound on the angle between the vectors).  

\begin{lemma}
Let $\bar{L}$ be the $(n-|S|)\times(n-|S|)$ Laplacian matrix for a connected network, and let  $\Delta = \hbox{diag}\left(\alpha_1, \alpha_2, \ldots, \alpha_{n-|S|}\right)$, where $0 \le \alpha_i \le |S|$ for all $1 \le i \le n-|S|$.  Let $\mathbf{x}$ be the positive eigenvector corresponding to the smallest eigenvalue of $L_g = \bar{L}+\Delta$, normalized so that $\|\mathbf{x}\|_{\infty} = 1$.  Then the smallest  element of $\mathbf{x}$ satisfies
\begin{equation}
x_{min} \geq 1 - \frac{2\sqrt{|S||\partial S|}}{\lambda_2(\bar{L})},
\label{eqn:eige12}
\end{equation}
where $|\partial S| \triangleq \sum_{i = 1}^{n-|S|}\alpha_i$.
 \label{lem:mainlem}
\end{lemma}

\begin{IEEEproof}
The eigenvector equation for $L_g$ is given by
\begin{equation}
L_g\mathbf{x} = (\bar{L}+\Delta)\mathbf{x}=\lambda \mathbf{x}.
\label{eqn:eige}
\end{equation}
Project the eigenvector $\mathbf{x}$ onto the subspace spanned by $\mathbf{v}_1$ to obtain $\mathbf{x}=\gamma \mathbf{1}+\mathbf{d}$, where $\mathbf{d}$ is orthogonal to $\mathbf{v}_1$ and $\gamma=\frac{\mathbf{1}^T\mathbf{x}}{n-|S|}$.   Thus we can write
\begin{equation}
\mathbf{d}=\sum_{i=2}^{n-|S|}\delta_i\mathbf{v}_i,
\label{eqn:d_decomp}
\end{equation}
for some real numbers $\delta_2, \delta_3, \ldots, \delta_{n-|S|}$.  Substituting $\mathbf{x}=\gamma \mathbf{1}+\mathbf{d}$ into \eqref{eqn:eige} and rearranging gives
\begin{equation}
\bar{L}\mathbf{d}=\underbrace{(\lambda I-\Delta)\mathbf{x}}_{\triangleq \mathbf{z}}.
\label{eqn:eige2}
\end{equation}
Multiplying  both sides of \eqref{eqn:eige2} by $\mathbf{1}^T$ yields $0 = \mathbf{1}^T\mathbf{z}$, and thus $\mathbf{z}$ is also orthogonal to $\mathbf{v}_1$.  Writing $\mathbf{z}=\sum_{i=2}^{n-|S|}\varphi_i\mathbf{v}_i$ for some constants $\varphi_2, \varphi_3, \ldots, \varphi_{n-|S|}$ and substituting this and \eqref{eqn:d_decomp} into \eqref{eqn:eige2}, we have
\begin{equation*}
\bar{L}\mathbf{d}=\sum_{i=2}^{n-|S|}\delta_i\bar{L}\mathbf{v}_i=\sum_{i=2}^{n-|S|}\delta_i\lambda_i(\bar{L})\mathbf{v}_i=\sum_{i=2}^{n-|S|}\varphi_i\mathbf{v}_i,
\end{equation*}
which gives $\delta_i=\frac{\varphi_i}{\lambda_i(\bar{L})}$ by the linear independence of the eigenvectors $\mathbf{v}_2, \ldots, \mathbf{v}_{n-|S|}$. Thus we can write $\mathbf{d}=\sum_{i=2}^{n-|S|} \frac{\varphi_i}{\lambda_i(\bar{L})}\mathbf{v}_i$ with $2$-norm given by
\begin{equation}
\|\mathbf{d}\|_2^2 = \sum_{i=2}^{n-|S|} \left( \frac{\varphi_i}{\lambda_i(\bar{L})} \right) ^2\leq \frac{1}{\lambda_2(\bar{L})^2}\sum_{i=2}^{n-|S|}\varphi_i^2=\frac{||\mathbf{z}||_2^2}{\lambda_2(\bar{L})^2}.
\label{eqn:eige5}
\end{equation}
From the definition of $\mathbf{z}$ in \eqref{eqn:eige2}, we have 
\begin{equation*}
\|\mathbf{z}\|_2^2 = \sum_{i=1}^{n-|S|}(\lambda - \alpha_i)^2 x_i^2 
\le \sum_{i=1}^{n-|S|}(\lambda - \alpha_i)^2 = (n-|S|)\lambda^2 - 2\lambda|\partial S|  + \sum_{i=1}^{n-|S|}\alpha_i^2.
\end{equation*}
Applying \eqref{eqn:maineq}, $|\partial S| \le |S|(n-|S|)$, and the fact that $\alpha_i \le |S|$ for all $1 \le i \le n-|S|$, we obtain
\begin{equation*}
\|\mathbf{z}\|_2^2 \leq (n-|S|)\frac{|\partial S|^2}{(n-|S|)^2} - 2\lambda|\partial S|  +  |S||\partial S| \leq 2|S||\partial S|.
\end{equation*}
Combining this with \eqref{eqn:eige5} yields
\begin{equation}
\|\mathbf{d}\|_2^2\leq \frac{2|S||\partial S|}{\lambda_2(\bar{L})^2}.
\label{eqn:eige8}
\end{equation}
Next, from $\mathbf{d}=\mathbf{x}-\gamma\mathbf{1}$ we have
\begin{equation}
\|\mathbf{d}\|_2^2\geq (x_{max}-\gamma)^2+(\gamma-x_{min})^2=(1-\gamma)^2+(\gamma-x_{min})^2.
\label{eqn:eige9}
\end{equation}
The right hand side of \eqref{eqn:eige9} achieves its minimum when  $\gamma=\frac{1+x_{min}}{2}$. Substituting this value and rearranging gives $x_{min}\geq 1-\sqrt{2}\|\mathbf{d}\|_2$, which yields the desired result when combined with \eqref{eqn:eige8}.
\end{IEEEproof}

The above result, in conjunction with Theorem~\ref{thm:main}, allows us to characterize graphs where the bounds in \eqref{eqn:maineq} become asymptotically tight. 

\begin{theorem}
Consider a sequence of connected graphs $\mathcal{G}_n$, $n \in \mathbb{Z}_{+}$, where $n$ indicates the number of nodes.  Consider an associated sequence of grounded nodes $S_n$, $n \in \mathbb{Z}_{+}$.  Let $\bar{L}_n$ denote the Laplacian matrix induced by the non-grounded nodes in each graph $\mathcal{G}_n$, and let $\lambda_n$ denote the smallest eigenvalue of the grounded Laplacian for the graph $\mathcal{G}_n$ with grounded set $S_n$.  Then:
\begin{enumerate}
\item If $\limsup_{n\rightarrow\infty}\frac{2\sqrt{|S_n||\partial S_n|}}{\lambda_2(\bar{L}_n)} < 1$, then $\lambda_n = \Theta\left(\frac{|\partial S_n|}{n-|S_n|}\right)$.
\item If $\lim_{n\rightarrow\infty}\frac{\sqrt{|S_n||\partial S_n|}}{\lambda_2(\bar{L}_n)} = 0$, then $(1-o(1))\frac{|\partial S_n|}{n-|S_n|} \le \lambda_n \le \frac{|\partial S_n|}{n-|S_n|}$.
\end{enumerate}
\label{thm:lambda_tight}
\end{theorem}

In the next sections, we will apply the results from this section to study the smallest eigenvalue of the grounded Laplacian of Erdos-Renyi and $d$-regular random graphs.

%%%%%%%%%%%%%%%%%%%%%%%%%%%%%%%%%%%%%%%%%%%%%%%%%%%%%%%%%%
%%%%%%%%%%%%%%%%%%%%%%%%%%%%%%%%%%%%%%%%%%%%%%%%%%%%%%%%%%%

\section{Analysis of Erdos-Renyi Random Graphs}

\begin{definition}
An Erdos-Renyi (ER) random graph, denoted $\mathcal{G}(n,p)$, is a graph on $n$ nodes where each possible edge between two distinct vertices is present independently with probability $p$ (which could be a function of $n$).  Equivalently, an ER random graph can be viewed as a probability space $(\Omega_n, \mathcal{F}_n, \mathbb{P}_n)$, where the sample space $\Omega_n$ consists of all possible graphs on $n$ nodes, the $\sigma$-algebra $\mathcal{F}_n$ is the power set of $\Omega_n$, and the probability measure $\mathbb{P}_n$ assigns a probability of $p^{|\mathcal{E}|}(1-p)^{\binom{n}{2}-|\mathcal{E}|}$ to each graph with $|\mathcal{E}|$ edges. 
\end{definition}

\begin{definition}
For an ER random graph, we say that a property holds {\it asymptotically almost surely} if the probability  of the set of graphs with that property (over the probability space $(\Omega_n, \mathcal{F}_n, \mathbb{P}_n)$) goes to $1$ as $n \rightarrow \infty$.  For a given graph function $f: \Omega_n \rightarrow \mathbb{R}_{\ge 0}$ and another function $g: \mathbb{N} \rightarrow \mathbb{R}_{\ge 0}$, we say $f(\mathcal{G}(n,p)) \le (1+o(1))g(n)$ asymptotically almost surely if there exists some function $h(n) \in o(1)$ such that $f(\mathcal{G}(n,p)) \le (1+h(n))g(n)$ with probability tending to $1$ as $n \rightarrow \infty$.  
\end{definition}

We start by showing the following bounds on the degrees and isoperimetric constants of such graphs; while there exist bounds on these quantities for specific forms of $p$ (e.g., \cite{Perez,Benjamini,Ganesh,Cooper2007cover}), they do not cover the full range of probability functions considered by the following lemma.   The proof of this result is provided in the Appendix.

%%%%%%%%%%%%%%%%%%%%%%%%%%%%%%%%%%%%%%%%%%%%%%%%%%%%%%%%%%%%%%%%%%%%%%

\begin{lemma}
Consider the Erdos-Renyi random graph $\mathcal{G}(n,p)$, where the edge probability $p$ satisfies 
$\limsup_{n \rightarrow \infty} \frac{\ln{n}}{np} < 1$.  Fix any $\epsilon \in (0, \frac{1}{2}]$.  There exists a positive constant $\alpha$ (that depends on $p$) such that the minimum degree $d_{min}$, maximum degree $d_{max}$ and isoperimetric constant $i(\mathcal{G})$ satisfy
$$
\alpha{np} \le i(\mathcal{G}) \le d_{min} \le d_{max} \le {np}\left(1 + \sqrt{3}\left(\frac{\ln n}{np}\right)^{\frac{1}{2}-\epsilon}\right).
$$
asymptotically almost surely.
\label{lem:degree_iso_bounds}
\end{lemma}

\begin{remark}
Note that the probability functions captured by the above lemma include the special cases where $p$ is a constant and where $p(n) = \frac{c\ln{n}}{n}$ for constant $c > 1$.  The above results generalize the bounds on the degrees and the isoperimetric constant in \cite{Perez,Benjamini,Cooper2007cover} where probability functions of the form $\frac{c\ln{n}}{n}$ were studied, although the bounding constants provided in those works will be generally tighter than the ones provided above due to the special case analysis. Further note that when $\ln{n} = o(np)$, the upper bound on the maximum degree becomes $np(1+o(1))$.
\end{remark}

Lemma~\ref{lem:degree_iso_bounds} and the lower bound \eqref{eqn:lower_bound_lambda_2_iso} immediately lead to the following corollary.

\begin{corollary}
Consider the Erdos-Renyi random graph $\mathcal{G}(n,p)$, where the edge probability $p$ satisfies
$\limsup_{n \rightarrow \infty} \frac{\ln{n}}{np} < 1$.  Then there exists a positive constant $\gamma$ (that depends on $p$) such that the algebraic connectivity $\lambda_2(\mathcal{G})$ satisfies $\lambda_2(L) \ge \gamma{np}$ asymptotically almost surely.
\label{cor:algebraic_connectivity_bound}
\end{corollary}

%%%%%%%%%%%%%%%%%%%%%%%%%%%%%%%%%%%%%%%%%%%%%%%%%%%%%%%%%%%%%%%%%%%%%%
With the above results in hand, we are now in place to prove the following fact about the smallest eigenvalue of the grounded Laplacian matrix for Erdos-Renyi random graphs.  We omit the dependence of $S$ and $\lambda$ on $n$ for notational convenience.

\begin{theorem}
Consider the Erdos-Renyi random graph $\mathcal{G}(n,p)$, where the edge probability $p$ satisfies
$\limsup_{n \rightarrow \infty} \frac{\ln{n}}{np} < 1$.  Let $S$ be a set of grounded nodes chosen uniformly at random with $|S| = o(\sqrt{np})$.  Then the smallest eigenvalue $\lambda$ of the grounded Laplacian satisfies
$$
(1-o(1))|S|p \le \lambda \le (1+o(1))|S|p
$$
asymptotically almost surely.
\label{thm:erdos1}
\end{theorem}

\begin{IEEEproof}
For probability functions satisfying the conditions in the theorem, Lemma~\ref{lem:degree_iso_bounds} indicates for any set $S$ of grounded nodes,  $|\partial S| \le |S|d_{max} \le \beta|S|np$ asymptotically almost surely for some positive constant $\beta$.  Let $\bar{L}$ be the Laplacian matrix for the graph induced by the non-grounded nodes (i.e., the graph obtained by removing all grounded nodes and their incident edges).  From \cite{Fiedler}, we have $\lambda_2(\bar{L}) \ge \lambda_2(L) - |S|$.  Combining this with Corollary~\ref{cor:algebraic_connectivity_bound}, we obtain
$$
\frac{\sqrt{|S||\partial S|}}{\lambda_2(\bar{L})} \le \frac{|S|\sqrt{\beta{np}}}{\gamma{np} - |S|} = o(1)
$$
asymptotically almost surely when $|S| = o(\sqrt{np})$.  From Lemma~\ref{lem:mainlem} and Theorem~\ref{thm:main}, we have $(1-o(1))\frac{|\partial S|}{n-|S|} \le \lambda \le \frac{|\partial S|}{n-|S|}$ asymptotically almost surely.

Next, consider the random variable $|\partial S|$; there are $|S|(n-|S|)$ possible edges between $S$ and $\mathcal{V}\setminus S$, each appearing independently with probability $p$, and thus $|\partial S|$ is a Binomial random variable with $|S|(n-|S|)$ trials.  For all $0<\alpha<1$ we have the concentration inequalities \cite{Mitzenmacher}
\begin{equation}
\begin{split}
\mathbf{Pr}(|\partial S| &\geq (1+\alpha)\mathbb{E}[|\partial S|])\leq e^{\frac{-\mathbb{E}[|\partial S|] \alpha^2}{3}}\\
\mathbf{Pr}(|\partial S| &\leq (1-\alpha)\mathbb{E}[|\partial S|])\leq e^{\frac{-\mathbb{E}[|\partial S|] \alpha^2}{2}}.
\end{split}
\label{eqn:2ineq}
\end{equation}
We know that $\mathbb{E}[|\partial S|]= |S|(n-|S|)p$. Consider $\alpha=\frac{1}{\sqrt[4]{\ln{n}}}$ which causes the upper bound in the first expression to become $\exp(-\frac{|S|(n-|S|)p}{3\sqrt{\ln{n}}})$.  Since $|S|(n-|S|)$ is lower bounded by $n-1$ and $np > \ln{n}$ for sufficiently large $n$, the bounds in \eqref{eqn:2ineq} asymptotically go to zero. Thus,
\begin{equation*}
(1 - o(1))|S|(n-|S|)p \le |\partial S| \le (1 + o(1))|S|(n-|S|)p,
\end{equation*}
asymptotically almost surely. Substituting into  the bounds for $\lambda$, we obtain the desired result.
\end{IEEEproof}

\section{Random $d$-Regular Graphs}
  \label{sec:random_d_regular}
  We now consider random $d$-regular graphs, defined as follows, and characterize the smallest eigenvalue of the grounded Laplacian for such graphs.

\begin{definition}
For any $n \in \mathbb{N}$, let $d \in \mathbb{N}$ be such that $3\leq d<n$ and $dn$ is an even number. Define $\Omega_{n,d}$ to be the set of all $d$-regular graphs on $n$ nodes. Define the probability space $(\Omega_{n,d}, \mathcal{F}_{n,d}, \mathbb{P}_{n,d})$, where the $\sigma$-algebra $\mathcal{F}_{n,d}$ is the power set of $\Omega_{n,d}$, and $\mathbb{P}_{n,d}$ is a probability measure assigning equal probability to every element of $\Omega_{n,d}$.  An element of $\Omega_{n,d}$ drawn according to $\mathbb{P}_{n,d}$ is called a random $d$-regular graph, and denoted by $\mathcal{G}_{n,d}$ \cite{Bollobas}.
\end{definition}
  
Let $\lambda'_1(A)\leq \lambda'_2(A)\leq ... \leq \lambda'_n(A)$ be the eigenvalues of the adjacency matrix of any given graph $\mathcal{G}$; note that $\lambda'_n(A)=d$ for $d$-regular graphs. Define $\lambda'(\mathcal{G})=\max \{|\lambda'_1(A)|,|\lambda'_{n-1}(A)|\}$.  It was shown in \cite{Friedman} that for any $\epsilon > 0$, 
\begin{equation}
  \lambda'(\mathcal{G}_{n,d})\leq 2\sqrt{d-1}+\epsilon
  \label{eqn:raman}
  \end{equation}
asymptotically almost surely.  As the Laplacian for the graph is given by $L = D - A = dI - A$, for any $\epsilon > 0$, the algebraic connectivity of a random $d$-regular graph satisfies
\begin{equation*}
\lambda_2(L) \geq d-2\sqrt{d-1}-\epsilon,
\end{equation*}
asymptotically almost surely.  
On the other hand we know that $\lambda_2(\bar{L})\geq \lambda_2(L)-|S|$ \cite{Fiedler}. Thus for a random $d$-regular graph with a single grounded node, we have
  \begin{equation*}
    \frac{2\sqrt{|S||\partial S|}}{\lambda_2(\bar{L})} =\frac{2\sqrt{d}}{\lambda_2(\bar{L})}\leq  \frac{2\sqrt{d}}{d-2\sqrt{d-1}-\epsilon-1}< 1
  \end{equation*}
  asymptotically almost surely for sufficiently large $d$ and sufficiently small $\epsilon$.   Lemma~\ref{lem:mainlem} and Theorem~\ref{thm:main} then yield the following result.
  
\begin{theorem}
Let $\mathcal{G}$ be a random $d$-regular graph on $n$ vertices with a single grounded node.  Then for sufficiently large $d$, the smallest eigenvalue of the grounded Laplacian satisfies $\lambda = \Theta\left(\frac{d}{n}\right)$ asymptotically almost surely.
\end{theorem}

%%%%%%%%%%%%%%%%%%%%%%%%%%%%%%%%%%%%%%%%%%%%%%%%%%%%%%%%%%%%%%%%%%%%%%%%%%%%%%%

\section{conclusion}
\label{sec:conclusion}
We studied the smallest eigenvalue of grounded Laplacian matrices, and provided bounds on this eigenvalue in terms of the number of edges between the grounded nodes and the rest of the network, bottlenecks in the network, and the smallest component of the eigenvector for the smallest eigenvalue.  We showed that our bounds are tight when the smallest eigenvector component is close to the largest component, and provided graph-theoretic conditions that cause the smallest component to converge to the largest component.  An outcome of our analysis is tight bounds for Erdos-Renyi random graphs and $d$-regular random graphs.  A rich avenue for future research is to extend and apply our results to other classes of random graphs.

\bibliographystyle{IEEEtran}
\bibliography{refs}

\begin{appendix}

\subsection{Proof of Lemma~\ref{lem:degree_iso_bounds}}
\begin{IEEEproof}
The degree bounds are readily obtained from classical concentration inequalities.  Specifically, let $d$ denote the degree of a given vertex.  Note that $d$ is a Binomial random variable with parameters $n-1$ and $p$, with expected value $\mathbb{E}[d] = (n-1)p$.  Now, for any $0 < \beta \le \sqrt{3}$ we have\footnote{The statement of this concentration inequality in \cite{Mitzenmacher} has $0 < \beta \le 1$, but the improved upper bound of $\sqrt{3}$ can be obtained from the same proof {\it mutatis mutandis}.} \cite{Mitzenmacher}
$$
\mathbf{Pr}(d \geq (1+\beta)\mathbb{E}[d])\leq e^{\frac{-\mathbb{E}[d] \beta^2}{3}}.
$$
Choose $\beta = \sqrt{3}\left(\frac{\ln{n}}{np}\right)^{\frac{1}{2}-\epsilon}$, which is at most $\sqrt{3}$ for probability functions satisfying the conditions in the lemma and for sufficiently large $n$.   Substituting into the above expression, we have 
\begin{equation*}
\mathbf{Pr}(d \geq (1+\beta)\mathbb{E}[d]) \leq e^{-(n-1)p \left(\frac{\ln{n}}{np}\right)^{1-2\epsilon}} =  O\left(e^{-\ln{n}\left(\frac{\ln{n}}{np}\right)^{-2\epsilon}}\right).
\end{equation*}
To show that the maximum degree is smaller than the given bound asymptotically almost surely, we show that all vertices have degree less than the given bound with probability tending to $1$.  By the union bound, the probability that at least one vertex has degree larger than $(1+\beta)\mathbb{E}[d]$ is upper bounded by
$$
n\mathbf{Pr}(d \geq (1+\beta)\mathbb{E}[d]) = O\left(e^{\ln{n}-\ln{n}\left(\frac{\ln{n}}{np}\right)^{-2\epsilon}}\right).
$$
Since $\limsup_{n\rightarrow\infty}\frac{\ln{n}}{np} < 1$, the above expression goes to zero as $n \rightarrow \infty$, proving the upper bound on the maximum degree.

We now show the lower bound for $i(\mathcal{G})$.  Specifically, we will show that for $p$ satisfying the properties in the lemma, almost every graph has the property that  all sets of vertices of size $s$, $1 \le s  \le \lfloor\frac{n}{2}\rfloor$, have at least $\alpha{snp}$ edges leaving that set, for some constant $\alpha$ that we will specify later.  For any specific set $\mathcal{S}$ of vertices of size $s$, the probability that $\mathcal{S}$ has $\lfloor\alpha{snp}\rfloor$ or fewer edges leaving the set is $\sum_{j = 0}^{\lfloor\alpha{snp}\rfloor}\binom{s(n-s)}{j}p^j(1-p)^{s(n-s)-j}$.  Let $E_s$ denote the event that at least one set of vertices of size $s$ has $\lfloor\alpha{snp}\rfloor$ or fewer edges leaving the set.  Then
\begin{equation}
\mathbf{Pr}\left[E_s\right] \le \binom{n}{s}\sum_{j = 0}^{\lfloor\alpha{snp}\rfloor}\binom{s(n-s)}{j}p^j(1-p)^{s(n-s)-j}.
\label{eqn:prob_Es}
\end{equation}
Note that for $1 \le j \le \lfloor\alpha{snp}\rfloor$,
\begin{align*}
\frac{\binom{s(n-s)}{j}p^j(1-p)^{s(n-s)-j}}{\binom{s(n-s)}{j-1}p^{j-1}(1-p)^{s(n-s)-j+1}} = \frac{s(n-s) - j + 1}{j}\frac{p}{1-p} &\ge \frac{s(n-s) -\alpha{s}{np}}{\alpha{s}{np}}\frac{p}{1-p} \\
&\ge \frac{1-2\alpha{p}}{2\alpha}\frac{1}{1-p} \ge \frac{1}{2\alpha}, 
\end{align*}
for $s \le \lfloor\frac{n}{2}\rfloor$ and $2\alpha < 1$ (which will be satisfied by our eventual choice of $\alpha$).  Thus, there exists some constant $r > 0$ such that 
\begin{equation*}
\sum_{j = 0}^{\lfloor\alpha{snp}\rfloor}\binom{s(n-s)}{j}p^j(1-p)^{s(n-s)-j} \le r\binom{s(n-s)}{\lfloor\alpha{snp}\rfloor}p^{\lfloor\alpha{snp}\rfloor}(1-p)^{s(n-s)-\lfloor\alpha{snp}\rfloor}.
\end{equation*}
Substituting into \eqref{eqn:prob_Es} and using the fact that $\binom{n}{k} \le \left(\frac{ne}{k}\right)^{k}$, we have
\begin{align}
\mathbf{Pr}\left[E_s\right] &\le r\left(\frac{ne}{s}\right)^{s}\left(\frac{s(n-s)ep}{\alpha{snp}}\right)^{\alpha{snp}}e^{-p(s(n-s)-\alpha{snp})} \nonumber \\
&\le re^{s\ln\frac{ne}{s}}\left(\frac{e}{\alpha}\right)^{\alpha{snp}}e^{-p(s(n-s)-\alpha{snp})} \nonumber \\
&= re^{sh(s)}, \label{eqn:prob_e_bound}
\end{align}
where
\begin{equation}
h(s) = 1 + np\underbrace{\left(\frac{\ln{n}}{np} + \alpha - \alpha\ln{\alpha} + \alpha{p} -1\right)}_{\Gamma(\alpha)} + ps - \ln{s}.
\label{eqn:hs}
\end{equation}
Noting that $h(s)$ is decreasing in $s$ until $s = \frac{1}{p}$ and increasing afterwards, we have 
\begin{align*}
h(s) &\le \max\left\{h(1), h\left(\frac{n}{2}\right)\right\} \\
&= \max\left\{1 +p + np\Gamma(\alpha), 1 + \ln{2} + np\left(\Gamma(\alpha) - \frac{\ln{n}}{np} + \frac{1}{2}\right)\right\}.
\end{align*}
From \eqref{eqn:hs},  $\Gamma(\alpha)$ is increasing in $\alpha$ for $0 \le \alpha < 1$, with $\Gamma(0) = \frac{\ln{n}}{np} - 1$ being negative and bounded away from $0$ for sufficiently large $n$ (by the assumption on $p$ from the statement of the lemma).   Thus, there exists some sufficiently small positive constant $\alpha$  such that $h(s) \le -\bar{\alpha}np$ for some constant $\bar{\alpha} > 0$ and for sufficiently large $n$.  Thus \eqref{eqn:prob_e_bound} becomes $\mathbf{Pr}\left[E_s\right] \le re^{-s\bar{\alpha}np}$ for sufficiently large $n$.  

By the union bound, the probability that $i(\mathcal{G}) < \alpha{np}$ is upper bounded by the sum of the probabilities of the events $E_s$ for $1 \le s \le \lfloor\frac{n}{2}\rfloor$.  Using the above expression, we have
$$
\sum_{s =1}^{\lfloor\frac{n}{2}\rfloor}\mathbf{Pr}[E_s] \le r\sum_{s = 1}^{\lfloor\frac{n}{2}\rfloor}e^{-s\bar{\alpha}np} \le r\sum_{s = 1}^{\infty}e^{-s\bar{\alpha}np} = r\frac{e^{-\bar{\alpha}np}}{1-e^{-\bar{\alpha}np}}
$$
which goes to $0$ as $n \rightarrow \infty$.  Thus, we have $i(\mathcal{G}) \ge \alpha{np}$ asymptotically almost surely.
\end{IEEEproof}

\end{appendix}

\end{document}